\documentclass[smallextended]{svjour3}
\usepackage{latexsym,amssymb,amsmath}
\smartqed
\usepackage{graphicx}

\RequirePackage{fix-cm}
\usepackage[section]{placeins}

\numberwithin{equation}{section}
\begin{document}
\title{\bf Piecewise Linear Vector Optimization Problems on Locally Convex Hausdorff Topological Vector Spaces}
\author{Nguyen Ngoc Luan} 
\institute{Nguyen Ngoc Luan \at Department of Mathematics and Informatics, Hanoi National University of Education, 136 Xuan Thuy, Hanoi, Vietnam\\
\email{luannn@hnue.edu.vn}}
\date{Received: date / Accepted: date}

\titlerunning{Piecewise Linear Vector Optimization Problems}

\maketitle

\begin{abstract} Piecewise linear vector optimization problems in a locally convex Hausdorff topological vector spaces setting are considered in this paper. The efficient solution set of these problems are shown to be the unions of finitely many semi-closed generalized polyhedral convex sets. If, in addition, the problem is convex, then the efficient solution set and the weakly efficient solution set are the unions of finitely many generalized polyhedral convex sets and they are connected by line segments. Our results develop the preceding ones of Zheng and Yang [Sci. China Ser. A. \textbf{51}, 1243--1256 (2008)], and Yang and Yen [J. Optim. Theory Appl. \textbf{147}, 113--124 (2010)], which were established in a normed spaces setting. 
	
\end{abstract}

\section{Introduction}
The intersection of a finite number of closed half-spaces of a finite-dimen\-sional Euclidean space is called a polyhedral convex set (a convex polyhedron in brief). Following Bonnans and Shapiro \cite[Definition~2.195]{Bonnans_Shapiro_2000}, we call a subset of a locally convex Hausdorff topological vector space (lcHtvs) a \textit{generalized polyhedral convex set} (or a \textit{generalized convex polyhedron}) if it is the intersection of finitely many closed half-spaces and a closed affine subspace of that topological vector space. If the affine subspace can be chosen as the whole space, the generalized polyhedral convex set (gpcs) is said to be a \textit{polyhedral convex set} (or a \textit{convex polyhedron}). 

\medskip
From now on, if not otherwise stated, $X$ and $Y$ are locally convex Hausdorff topological vector spaces. Similarly as in \cite{ZhengYang_2008}, we say that a mapping $f:X\to Y$ is a {\it piecewise linear function} (or a piecewise affine function) if there exist polyhedral convex sets $P_1, \dots, P_m$  in $X$, continuous linear mappings $T_1,\dots, T_m$  from $X$ to $Y$, and vectors  $b_1, \dots, b_m$ in $Y$ such that  $X=\bigcup\limits_{k=1}^m P_k$ and $f(x)=T_kx + b_k$ for all $x \in P_k$, $k=1,\dots,m$. 

\medskip
We call a subset $K \subset Y$ a cone if $tK \subset K$ for all $t>0$. 

\medskip 
Given a piecewise linear function $f: X \rightarrow Y$, a generalized polyhedral convex set $D \subset X$, and a polyhedral convex cone $K \subset Y$ with $K\neq Y$, we consider the {\it piecewise linear vector optimization problem}   
\begin{equation*}
{\rm (VP)} \qquad \quad {\rm Min}_K \big\{f(x) \mid x \in D\big\}.
\end{equation*}
In the terminology of \cite[p.~341]{Giannessi_1984}, one says that $f$ is a $K$-\textit{function}\footnote{Some authors use the term ``\textit{$K$-convex function}" instead of $K$-function.} on $D$ if 
\begin{equation*}
(1-\lambda)f(x_1) + \lambda f(x_2) - f ((1-\lambda) x_1 + \lambda x_2) \in K
\end{equation*} 
for any $x_1, x_2$ in $D$ and $\lambda \in [0,1]$. It is clear that if $f$ is linear, then it is a $K$-function on $D$. In the case where $Y=\mathbb{R}$ and $K=\mathbb{R}_+$, $f$ is a $K$-function on~$D$ if and only if $f$ is convex on $D$. When $f$ is a $K$-function on~$D$, we say that {\rm (VP)} is a {\it convex problem}. 

\medskip
A point $y$ from a subset $Q\subset Y$ is called an \textit{efficient point} (resp., a \textit{weakly efficient point}) of $Q$ if there is no $y' \in Q$ such that $y - y' \in K\setminus{\ell(K)}$ (resp., $y-y' \in {\rm int}\,K$), where $\ell(K):=K \cap (-K)$. The efficient point set and the weakly efficient point set of $Q$ are denoted respectively by $E(Q|K)$ and $E^w(Q|K)$. Clearly, in the case where $K$ is a pointed cone, i.e., $\ell(K)=\left\{ 0 \right\}$, a point $y \in Q$ belongs to $E(Q|K)$ if and only if $y-y' \notin K\setminus\{0\}$ for every $y' \in Q$. Since ${\rm int}\,K \subset K\setminus{\ell(K)}$ by \cite[Proposition~2.2]{Luan_2016}, one gets $E(Q|K) \subset E^w(Q|K)$. 

\medskip
A vector $u \in D$ is called an \textit{efficient solution} (resp.,  a \textit{weakly efficient solution}) of {\rm (VP)}  if $f(u) \in E(f(D)|K)$ (resp., $f(u) \in E^w(f(D)|K)$). The efficient solution set and the weakly efficient solution set of {\rm (VP)} are denoted respectively by ${\rm Sol(VP)}$ and ${\rm Sol}^w{\rm (VP)}$. Since $E(Q|K) \subset E^w(Q|K)$, one has ${\rm Sol(VP)} \subset {\rm Sol}^w{\rm (VP)}$.  

\medskip
The study of the structures and characteristic properties of these solution sets is useful in the design of efficient algorithms for solving {\rm (VP)}. 

\medskip
Zheng and Yang \cite{ZhengYang_2008} have proved that if $X,\, Y$ are normed spaces and $D$ is a polyhedral convex set, then ${\rm Sol}^w{\rm (VP)}$ is the union of finitely many polyhedral convex sets. If $f$ is a $K$-function, then ${\rm Sol}^w{\rm (VP)}$ is connected by line segments. 

\medskip
In order to describe the structure of ${\rm Sol(VP)}$ and obtain sufficient conditions for its connectedness, Yang and Yen \cite{Yang_Yen_2010} have applied the {\it image space approach} \cite{Giannessi_1984,Giannessi_1987} to optimization problems and variational systems and proposed the notion of {\it semi-closed polyhedral convex set}. On account of \cite[Theorem~2.1]{Yang_Yen_2010}, if $X, \, Y$ are normed spaces, $Y$ is of finite dimension, $K \subset Y$ is a pointed cone, and $D \subset X$ is a polyhedral convex set, then ${\rm Sol(VP)}$ is the union of finitely many semi-closed polyhedra. In this setting, if $f$ is a $K$-function on $D$, then ${\rm Sol(VP)}$ is the union of finitely many polyhedra and it is connected by line segments; see \cite[Theorem~2.2]{Yang_Yen_2010}. Observe that the main tool for proving the latter results is the representation formula for convex polyhedra in $\mathbb{R}^n$ via a finite number of points and a finite number of  directions (see, e.g., \cite[Theorem~2.12]{Klee_1959} and \cite[Theorem~19.1]{Rockafellar_1970}). This celebrated result is attributed \cite[p. 427]{Rockafellar_1970} primarily to Minkowski \cite{Minkowski_1910} and Weyl \cite{Weyl_1935,Weyl_1953}. 

\medskip
Theorem~2.3 of \cite{Yang_Yen_2010} is an infinite-dimensional version of the classic Arrow, Barankin and Blackwell theorem (the ABB theorem; see \cite{ABB_1953,Luc_1989,Luc_2016}), which says that the efficient solution set of a finite-dimensional linear vector optimization problem is the union of finitely many polyhedral convex sets and it is connected by line segments.

\medskip
In \cite{Fang_Meng_Yang_2012}, Fang, Meng, and Yang studied multiobjective optimization problems with either continuous or discontinuous piecewise linear objective functions and polyhedral convex constraint sets. They obtained an algebraic representation of a semi-closed polyhedron and apply it to show that the image of a semi-closed polyhedron under a continuous linear function is always a semi-closed polyhedron. They proposed an algorithm for finding the Pareto point set of a continuous piecewise linear bi-criteria program and generalized it to the discontinuous case. The authors applied that algorithm to solve discontinuous bi-criteria portfolio selection problems with an $\ell_\infty$ risk measure and transaction costs. Some examples with the historical data of the Hong Kong Stock Exchange are discussed. Other results in this direction were given in~\cite{Fang_Huang_Yang_2012} and \cite{Fang_Meng_Yang_2015}. Later, Zheng and Ng \cite{Zheng_Ng_2014} have studied the metric subregularity of piecewise polyhedral multifunctions and applied this property to piecewise linear multiobjective optimization. 

\medskip
Very recently, in a locally convex Hausdorff topological vector spaces setting, by the use of a representation formula for generalized polyhedral convex sets, Luan and Yen~\cite{Luan_Yen_2015} have obtained solution existence theorems for generalized linear programming problems, a scalarization formula for the weakly efficient solution set of a generalized linear vector optimization problem, and proved that the latter is the union of finitely many generalized polyhedral convex sets. In \cite{Luan_2016}, where the relative interior of the dual cone of a polyhedral convex cone is described, a similar result is given for the corresponding efficient solution set. Moreover, it is shown that both efficient solution set and weakly efficient solution set of a generalized linear vector optimization problem are connected by line segments. Thus, the ABB theorem in linear vector optimization has been extended to the lcHtvs setting. 

\medskip
Various generalized polyhedral convex constructions in locally convex Hausdorff topological vector spaces can be found in the new paper of Luan, Yao, and Yen \cite{Luan_Yao_Yen_2016}. 

\medskip
It is well known that any infinite-dimensional normed space equipped with the {\it weak topology} is not metrizable, but it is a locally convex Hausdorff topological vector space. Similarly, the dual space of any infinite-dimensional normed space equipped with the {\it weak$^*$-topology} is not metrizable, but it is a locally convex Hausdorff topological vector space. Actually, the just mentioned two models provide us with the most typical examples of locally convex Hausdorff topological vector space, whose topologies cannot be given by norms. 

\medskip
The aim of the present paper is to revisit the results of Zheng and Yang~\cite{ZhengYang_2008}, Yang and Yen~\cite{Yang_Yen_2010} in a broader setting. Namely, instead of normed spaces, we will consider locally convex Hausdorff topological vector spaces. Apart from using some ideas and schemes of \cite{Yang_Yen_2010}, our investigation is based on the results of \cite{Luan_2016,Luan_Yao_Yen_2016,Luan_Yen_2015}.

\medskip
The paper organization is as follows. The next section gives some preliminaries. In Section 3, we prove that if the problem {\rm (VP)} is convex, then the solution sets are the unions of finitely many generalized polyhedral convex sets and they are connected by line segments. In Section 4, without imposing the convexity assumption on {\rm (VP)}, we show that ${\rm Sol}{\rm (VP)}$ is the union of finitely many semi-closed generalized polyhedral convex sets, while ${\rm Sol}^w{\rm (VP)}$ is the union of finitely many generalized polyhedral convex sets. Two illustrative examples are given in Sections 3 and 4.

\section{Preliminaries}
We now recall some concepts and results on generalized polyhedral convex sets and polyhedral convex cones. Let $X$ be a {\it locally convex Hausdorff topological vector space}. Denote by $X^*$ the dual space of $X$ and by $\langle x^*, x \rangle$ the value of~$x^* \in X^*$ at $x \in X$. 	

\begin{definition}\label{Def_gpcs} (See \cite[p.~133]{Bonnans_Shapiro_2000}) A subset $D \subset X$ is said to be a \textit{generalized polyhedral convex set}, or a \textit{generalized convex polyhedron}, if there exist $x^*_i \in X^*$, $\alpha_i \in \mathbb R$, $i=~1,2,\dots,p$, and a closed affine subspace $L \subset X$, such that 
\begin{equation}\label{eq_def_gpcs}
		D=\big\{ x \in X \mid x \in L,\ \langle x^*_i, x \rangle \leq \alpha_i,\  i=1,\dots,p\big\}.
\end{equation} 
If $D$ can be represented in the form \eqref{eq_def_gpcs} with $L=X$, then we say that it is a \textit{polyhedral convex set} (pcs), or a \textit{convex polyhedron}.
\end{definition}

If $L$ is a closed affine subspace of $X$ then, by \cite[Remark~2.196]{Bonnans_Shapiro_2000}, there exists a continuous surjective linear mapping $A$ from $X$ to a lcHtvs $Z$ and a vector $z \in Z$ such that $L=\big\{x \in X \mid Ax=z  \big\}$. So, one can rewrite \eqref{eq_def_gpcs} as follows
\begin{equation}\label{eq_def_gpcs_2}
D=\big\{ x \in X \mid Ax=z,\  \langle x^*_i, x \rangle \leq \alpha_i,\  i=1,\dots,p\big\}.
\end{equation}   

\medskip
Next lemmas, which were obtained in \cite{Luan_2016,Luan_Yao_Yen_2016,Luan_Yen_2015}, will be useful for our subsequent investigations.

\begin{lemma}\label{rep_gpcs_LY_2015}{\rm (See \cite[Theorem~2.7]{Luan_Yen_2015} and \cite[Lemma~2.12]{Luan_Yao_Yen_2016})} Suppose that $D$ is a nonempty subset of $X$. The set $D$ is generalized polyhedral convex if and only if there exist $u_1, \dots, u_k \in X$, $v_1, \dots, v_{\ell} \in X$, and a closed linear subspace  $X_0 \subset X$ such that
	\begin{equation}\label{rep_gpcs}
	\begin{aligned}
	D\!=\!\Bigg\{ \sum\limits_{i=1}^k \lambda_i u_i + \sum\limits_{j=1}^\ell \mu_j v_j \mid & \lambda_i \geq 0, \, \forall i=1,\dots,k,    \\ 
	&\sum\limits_{i=1}^k \lambda_i=1, \, \mu_j \geq 0,\, \forall j=1,\dots,\ell \Bigg\}+X_0.&
	\end{aligned}
	\end{equation}
In particular, $D$ is polyhedral convex if and only if $D$ admits a representation of the form \eqref{rep_gpcs}, where $X_0$ is a closed linear subspace of finite codimension.	
\end{lemma}

\begin{lemma}\label{sum_of_gpcs_LYY_2016}{\rm (See \cite[Propositions~2.10--2.11]{Luan_Yao_Yen_2016})} Suppose that $D_1, D_2$ are generalized polyhedral convex sets of $X$. If ${\rm aff}D_1$ is finite-dimensional or $D_1$ is polyhedral convex, then $D_1+D_2$ is a generalized polyhedral convex set.
\end{lemma}

\begin{lemma}\label{closed_sum_two_subspaces_LYY_2016}{\rm (See \cite[Lemma~2.13]{Luan_Yao_Yen_2016})} If $X_1$ and $X_2$ are linear subspaces of $X$ with $X_1$ being closed and finite-codimensional, then $X_1+X_2$ is closed and finite-codimensional.
\end{lemma}

\begin{lemma}\label{union_of_gpcs_LYY_2016}{\rm (See \cite[Corollary~2.16]{Luan_Yao_Yen_2016})} If a convex subset $D \subset X$ is the union of a finite number of generalized polyhedral convex sets (resp., of polyhedral convex sets) in $X$, then $D$ is generalized polyhedral convex (resp., polyhedral convex).   
\end{lemma}

\begin{lemma}\label{image_gpcs_LYY_2016}{\rm (See \cite[Proposition~2.7]{Luan_Yao_Yen_2016})} Let $T: X \rightarrow Y$ be a continuous linear mapping between locally convex Hausdorff topological vector spaces. If $Q \subset Y$ is a nonempty generalized polyhedral convex set, then $T^{-1}(Q)$ is a generalized polyhedral convex set.	
\end{lemma}

\begin{lemma}\label{cont_proj_LY_2015}{\rm (See \cite[Lemma~2.5]{Luan_Yen_2015})} If $W$ and $Z$ are Hausdorff finite-dimensional topological vector spaces of dimension $n$ and if $g: W \rightarrow Z$ is a linear bijective mapping, then $g$ is a homeomorphism.
\end{lemma}

\smallskip
According to Yang and Yen \cite{Yang_Yen_2010}, a subset of a normed space is called a semi-closed polyhedron if it is the intersection of a finite family of (closed or open) half-spaces. The following definition appears naturally in that spirit.
\begin{definition}\label{Def_smc_gpcs} A subset $D \subset X$ is said to be a \textit{semi-closed generalized polyhedral convex set}, or a \textit{semi-closed generalized convex polyhedron}, if there exist $x^*_i \in~X^*$, $\alpha_i \in \mathbb R$, $i=~1,2,\dots,q$, with a positive integer $p \leq q$, and a closed affine subspace $L \subset X$, such that 
\begin{equation}\label{eq_def_smc_gpcs}
	D\!=\!\big\{ x \in L \mid  \langle x^*_i, x \rangle \leq \alpha_i,\,  i=1,\dots,p;\, \langle x^*_i, x \rangle < \alpha_i,\,  i=p+1,\dots,q\big\}.
\end{equation}
If $D$ can be represented in the form \eqref{eq_def_smc_gpcs} with $L=X$, then we say that it is a \textit{semi-closed polyhedral convex set}, or a \textit{semi-closed convex polyhedron}. 
\end{definition}

\begin{remark} From Definitions~\ref{Def_gpcs} and \ref{Def_smc_gpcs} it follows that every polyhedral convex set is a semi-closed polyhedral convex set. It is not difficult to show that, if $D_1$ and $D_2$ are semi-closed polyhedral convex sets, then $D_1\setminus D_2$ is the union of finitely many semi-closed polyhedral convex sets.
\end{remark}

Suppose that $Y$ is a locally convex Hausdorff topological vector space and $K \subset Y$ is a polyhedral convex cone defined by
\begin{equation}\label{reps_K}
K=\big\{ y \in Y \mid \langle y^*_j, y \rangle \leq 0,\ j=1,\dots,q  \big\}, 
\end{equation}
where $y^*_j \in Y^*\setminus\{0\}$ for $j=1,\dots,q$. Note that
\begin{equation*}
\ell(K)=Y_0:=\big\{ y \in Y \mid \langle y^*_j, y \rangle = 0,\ j=1,\dots,q  \big\}
\end{equation*}
is a closed linear subspace of $Y$ and the codimension of $Y_0$, denoted by ${\rm codim}Y_0$, is finite. Moreover, there exists a finite-dimensional linear subspace $Y_1$ of $Y$, such that $Y=Y_0 + Y_1$ and $Y_0 \cap Y_1=\{0\} $. According to \cite[Theorem~1.21(b)]{Rudin_1991}, $Y_1$ is closed. Clearly, $K_1:=\big\{ y \in Y_1 \mid \langle y^*_j, y \rangle \leq 0, \; j=1,\dots,q\big\}$ is a pointed polyhedral convex cone in $Y_1$ and 
\begin{equation}\label{rep_for_K} 
K=Y_0+K_1.
\end{equation}

\medskip
The assertions of two following lemmas have been proved in \cite{Luan_2016}. 

\begin{lemma}\label{decomp_K_diff_lK} {\rm (See \cite[Lemma~2.1]{Luan_2016})} It holds that $K \setminus \ell(K)=Y_0+K_1 \setminus \{0\}$. 
\end{lemma}

\begin{lemma}\label{reps_intK}{\rm (See \cite[Proposition~2.2]{Luan_2016})} Let $K \subset Y$ be a polyhedral convex cone of the form \eqref{reps_K}. The following are valid:
\begin{description}
	\item{\rm (a)} The interior of $K$ has the representation 
	\begin{equation*}\label{reps_intK_eq}
	{\rm int}\,K=\big\{ y \in Y \mid \langle y^*_j, y \rangle < 0,\ j=1,\dots,q  \big\}.
	\end{equation*}
	\item{\rm (b)} The set $K \setminus{\ell(K)}$ is a convex cone and
	\begin{equation*}\label{reps_1b}
	K \setminus{\ell(K)}=\Big\{ y \in K \mid \text{there exists $j \in \{1,\dots,q\}$ such that } \langle y^*_j, y \rangle< 0  \Big\}.
	\end{equation*}	
\end{description}
\end{lemma}

The invariance of $K\setminus{\ell(K)}$ and ${\rm int}\,K$ w.r.t. a translation by a vector from~$K$ is described by the forthcoming lemma, which can be proved similarly as \cite[Proposition~4.3, p. 19]{Luc_1989} (see also \cite[Lemma~1.2 (iii)]{ChenHuangYang_2005}).  
\begin{lemma}\label{sum_cones} We have 
\begin{description}
	\item{\rm (a)} $K\setminus{\ell(K)}+K\subset K\setminus{\ell(K)}$;
	\item{\rm (b)} ${\rm int}\,K + K \subset {\rm int}\,K$.
\end{description} 	
\end{lemma}	

\section{Structure of the Solution Sets in the Convex Case}
One says that a subset $A \subset X$ is \textit{connected by line segments} if for any points $u, u'$  in~$A$, there are some points $u_1,\dots, u_r$ in $A$ with $u_1=u$ and $u_r=u'$ such that $[u_i, u_{i+1}] \subset A$ for $i=1,2,\dots,r-1$.
 
\medskip
The next result is an extension of \cite[Theorem~2.2]{Yang_Yen_2010} and \cite[Theorems~3.2 and 3.3]{ZhengYang_2008} to a locally convex Hausdorff topological vector spaces setting.

\begin{theorem}\label{struc_sol_convex} If $f$ is a $K$-function on $D$, the efficient solution set and the weakly efficient solution set of {\rm (VP)} are the unions of finitely many generalized polyhedral convex sets and they are connected by line segments.
\end{theorem}

The proof of this result is based on the next lemma.
\begin{lemma}\label{structure_solution_set_VLP_L_2016}{\rm (See \cite[Theorems~3.2--3.4]{Luan_2016} and \cite[Theorem~4.5]{Luan_Yen_2015})} If the function $f$ is linear, then ${\rm Sol(VP)}$ (resp., ${\rm Sol}^w{\rm (VP)}$) is the union of finitely many generalized polyhedral convex sets and it is connected by line segments.
\end{lemma}
\noindent
{\it Proof of Theorem~\ref{struc_sol_convex}} \ Without loss of generality, we may assume that $D$ is nonempty.

\medskip
\noindent
{\sc Claim 1.} \textit{The sum $f(D)+K$ is a polyhedral convex set.}

\smallskip
Following \cite[p.~18]{Luc_1989}, we define the epigraph of $f$ by the formula 
$${\rm epi}\,f=\{(x, y) \in X \times Y \mid y \in f(x)+K\}.$$
Since $f$ is a $K$-function, ${\rm epi}\,f$ is convex by \cite[Proposition~6.2, p.~29]{Luc_1989}. As $f(D)+K$ is the projection of ${\rm epi}\,f \cap (D \times Y)$ on $Y$, it is convex.  

For $k=1,\dots,m$, set $M_k:=f(D \cap P_k)$ and note that $$M_k=T_k(D \cap P_k) + b_k.$$ We will show that $M_k+K$, $k=1,\dots,m$, are pcs. Clearly, $D \cap P_k$ is a gpcs in~$X$. If $D \cap P_k$ is empty, then $M_k=\emptyset$; hence $M_k+K=\emptyset$ is a special polyhedral convex set. In the case where $D \cap P_k$ is nonempty, by the representation for gpcs (see Lemma \ref{rep_gpcs_LY_2015}) one can find $u_{k,1},\dots,u_{k,r_k}$, $v_{k,1},\dots,v_{k,s_k}$ in $X$ and a closed linear subspace  $X_{0,k} \subset X$ such that
	\begin{equation*}\label{rep_DcapPk}
	\begin{aligned}
	D\cap P_k\!=\!\Big\{ \sum\limits_{i=1}^{r_k}\lambda_i u_{k,i} + \sum\limits_{j=1}^{s_k} \mu_j v_{k,j} \mid &\, \lambda_i \geq 0,\, i=1,\dots,r_k,    \\ 
	&\, \sum\limits_{i=1}^{r_k} \lambda_i=1,\, \mu_j \geq 0,\, j=1,\dots,s_k  \Big\}\!+\!X_{0,k}.&
	\end{aligned}
	\end{equation*}
It is not difficult to show that $M_k=M'_k+T_k(X_{0,k})$, where
	\begin{equation*}\label{}
	\begin{aligned}
	M'_k\!=\!\Big\{ \sum\limits_{i=1}^{r_k} \lambda_i (Tu_{k,i}+b_k) + \sum\limits_{j=1}^{s_k} \mu_j Tv_{k,j} \mid &\, \lambda_i \geq 0,\, i=1,\dots,r_k,    \\ 
	&\sum\limits_{i=1}^{r_k}\lambda_i=1,\, \mu_j \geq 0,\, j=1,\dots,s_k \Big\}.&
	\end{aligned}
	\end{equation*}	
Combining this with the equality \eqref{rep_for_K}, one has $$M_k+K=M'_k+T_k(X_{0,k})+Y_0+K_1.$$ Since $T_k(X_{0,k})$ is a linear subspace of $Y$ and $Y_0$ is a closed finite-codimensional linear subspace of~$Y$, Lemma~\ref{closed_sum_two_subspaces_LYY_2016} shows that $Y_0+T_k(X_{0,k})$ is a closed subspace of $Y$ and ${\rm codim}\big(Y_0+T_k(X_{0,k})\big)$ is finite. Hence, by Lemma~\ref{rep_gpcs_LY_2015}, $M'_k+\big(Y_0+T_k(X_{0,k})\big)$ is a pcs in $Y$. As $K_1$ is a gpcs in $Y$, Lemma~\ref{sum_of_gpcs_LYY_2016} yields $K_1+ M'_k+Y_0+T_k(X_{0,k})$ is a polyhedral convex set. So, $M_k+K$ is a pcs in $Y$. 

In accordance with Lemma \ref{union_of_gpcs_LYY_2016}, since the convex set $f(D)+K$ is the union of  polyhedral convex sets $M_1+K, \dots, M_m+K$, we may conclude that $f(D)+K$ is a polyhedral convex set.   

\medskip
\noindent
{\sc Claim 2.} \textit{The sets $E(M+K|K)$ and $E^w(M+K|K)$, where $M:=f(D)$, are the unions of finitely many generalized polyhedral convex sets and they are connected by line segments.}

\smallskip
It suffices to apply Lemma~\ref{structure_solution_set_VLP_L_2016} to the problem  
\begin{equation*}
{\rm (VLP)} \qquad \quad {{\rm Min}_K} \left\{Gy \mid y \in M+K\right\}
\end{equation*}
with $G(y):=y$ for all $y \in Y$.

\medskip
\noindent
{\sc Claim 3.} \textit{A vector $u \in D$ belongs to {\rm Sol(VP)} (resp., belongs to ${\rm Sol}^w{\rm (VP)}$) if and only if $f(u) \in E(M+K|K)$ (resp., $f(u) \in E^w(M+K|K)$).}

\smallskip
Arguing similarly as in the proof of \cite[Proposition~7.10]{Luc_Ratiu_2014}, we can show that 
\begin{equation}\label{relation_1}
E(M|K) = M \cap E(M + K|K),\ E^w(M|K) = M \cap E^w(M + K|K).
\end{equation}
 For any $u \in D$, one has $f(u) \in M$. By definition, $u$ belongs to ${\rm Sol(VP)}$ (resp., to ${\rm Sol}^w{\rm (VP)}$) if and only if $f(u) \in E(M|K)$ (resp., $f(u) \in E^w(M|K)$). Hence, the assertions follow from \eqref{relation_1}. 

\medskip
\noindent
{\sc Claim 4.} \textit{The efficient solution set and the weakly efficient solution set of~{\rm (VP)} are the unions of finitely many generalized polyhedral convex sets.}

\smallskip
According to Claim~2, one can find generalized polyhedral convex sets $Q_1, \dots, Q_d$ in $Y$ such that $E(M+K|K)=\bigcup\limits\limits_{j=1}^d Q_j$. For each $k\in \{1,\dots,m\}$, by using Claim~3, one has
\begin{equation*}
\begin{aligned}
{\rm Sol(VP)} \cap P_k&=\big\{u \in D \cap P_k \mid f(u) \in E(M|K)\big\}\\
&=\big\{u \in D \cap P_k \mid f(u) \in E(M+K|K)\big\}\\
&=\bigcup\limits\limits_{j=1}^d  \big\{u \in D \cap P_k \mid f(u) \in Q_j\big\}\\
&=\bigcup\limits\limits_{j=1}^d  \big\{u \in D \cap P_k \mid  T_ku+b_k \in Q_j\big\}\\
&=\bigcup\limits\limits_{j=1}^d \Big( D \cap P_k \cap T_k^{-1}(-b_k+Q_j) \Big).
\end{aligned}
\end{equation*}
Since ${\rm Sol(VP)}=\bigcup\limits_{k=1}^m \big({\rm Sol(VP)} \cap P_k\big)$, we obtain 
\begin{equation}\label{union_for_Sol}
{\rm Sol(VP)}=\bigcup\limits_{k=1}^m \bigcup\limits\limits_{j=1}^d \Big( D \cap P_k \cap T_k^{-1}(-b_k+Q_j) \Big).
\end{equation}
For any $k \in \{1,\dots,m\}$ and $j \in \{1,\dots,d\}$, the set $-b_k+Q_j$ is generalized polyhedral convex by Lemma~\ref{sum_of_gpcs_LYY_2016}. Hence, since $T_k$ is a continuous linear mapping, by Lemma~\ref{image_gpcs_LYY_2016} we can assert that $T_k^{-1}(-b_k+Q_j)$ is a gpcs. This implies that $D \cap P_k \cap T_k^{-1}(-b_k+Q_j)$ is a gpcs. Then, formula~\eqref{union_for_Sol} justifies the fact that ${\rm Sol(VP)}$ is the union of finitely many generalized polyhedral convex sets. The assertion concerning ${\rm Sol}^w{\rm (VP)}$ can be proved similarly.

\medskip
\noindent
{\sc Claim 5.} \textit{The efficient solution set and the weakly efficient solution set of~{\rm (VP)} are connected by line segments.}

\smallskip
By analogy, it suffices to prove the assertion about the efficient solution set. Take any $u, u'$ in ${\rm Sol(VP)}$. By Claim~3, $f(u)$ and $f(u')$ are contained in $E(M+K|K)$. Since $E(M+K|K)$ is connected by line segments (see Claim~2), there exist $y_1,\dots,y_r$ in $Y$ such that $y_1=f(u)$, $y_r=f(u')$ and $$[y_i, y_{i+1}] \subset E(M+K|K) \quad (i=1,\dots,r-1).$$ Let $u_i \in D$ and $w_i \in K$ be such that $y_i=f(u_i)+w_i$ for $i=2,\dots,r-1$. Setting $u_1=u$, $w_1=0$, $u_r=u'$, and $w_r=0$, we have $y_i=f(u_i)+w_i$ for $i=1,\dots,r$.  Let $i \in \{1,\dots,r-1\}$ be chosen arbitrarily. To show that $[u_i, u_{i+1}] \subset {\rm Sol(VP)}$, we suppose the contrary: There is a $\lambda \in [0,1]$ such that the vector $u_\lambda:=(1-\lambda)u_i + \lambda u_{i+1}$ does not belong to ${\rm Sol(VP)}$. Then one can find $\bar x \in D$ and $\bar w \in K\setminus{\ell(K)}$ satisfying $f(u_\lambda)=f(\bar x)+\bar w$. As $f$ is a $K$-function on~$D$, there is $w \in K$ with
\begin{equation}\label{equality_1}
(1-\lambda) f(u_i) +\lambda f( u_{i+1} )=f((1-\lambda)u_i + \lambda u_{i+1})+w.
\end{equation}
Set $y_\lambda:=(1-\lambda) y_i +\lambda y_{i+1}$. On one hand, the inclusion $[y_i, y_{i+1}] \subset E(M+K|K)$ gives $y_\lambda\in E(M+K|K)$. On the other hand, 
\begin{equation}\label{equality_2}
\begin{aligned}
y_\lambda=&\, (1-\lambda) (f(u_i)+w_i) +\lambda (f(u_{i+1})+w_{i+1})\\
=&\, [(1-\lambda) f(u_i)+\lambda f(u_{i+1})] + (1-\lambda) w_i +\lambda w_{i+1}.
\end{aligned}
\end{equation}
Combining \eqref{equality_1}, \eqref{equality_2} with the equality $f(u_\lambda)=f(\bar x)+\bar w$, we obtain
\begin{equation*}
\begin{aligned}
y_\lambda &= [f((1-\lambda)u_i + \lambda u_{i+1})+w]+(1-\lambda) w_i +\lambda w_{i+1}.\\
&=[f(u_\lambda)+w] +(1-\lambda) w_i +\lambda w_{i+1}\\
&=f(\bar x)+\bar w + (w +(1-\lambda) w_i +\lambda w_{i+1}).
\end{aligned}
\end{equation*}
By the convexity of the cone $K$, one has $w +(1-\lambda) w_i +\lambda w_{i+1} \in K$. Therefore, as $\bar w \in K\setminus{\ell(K)}$, Lemma~\ref{sum_cones} shows that $\bar w + (w +(1-\lambda) w_i +\lambda w_{i+1}) \in K\setminus{\ell(K)}$. It follows that $y_\lambda - f(\bar x) \in K\setminus{\ell(K)}$. Hence, $y_\lambda\notin E(M+K|K)$, a contradiction. We have proved that the line segments $[u_i,u_{i+1}], i=1,\dots,r-1,$ which connect $u$ with $u'$, lie wholly in {\rm Sol(VP)}. 

The proof is complete.$\hfill\Box$

\begin{remark} In the proof of Theorem~\ref{struc_sol_convex}, we have used some ideas of Yang and Yen~\cite{Yang_Yen_2010}. For the case where $X, Y$ are normed spaces, $Y$ is finite-dimensional, $K \subset Y$ is a pointed cone, and $D$ is a polyhedral convex set, the assertions of Theorem~\ref{struc_sol_convex} about ${\rm Sol(VP)}$ recover \cite[ Theorem~2.2]{Yang_Yen_2010}. For the case where $X, Y$ are normed spaces and $D$ is a pcs, the assertions of Theorem~\ref{struc_sol_convex} about ${\rm Sol}^w{\rm (VP)}$ have been established by Zheng and Yang \cite[Theorems~3.2 and~3.3]{ZhengYang_2008}. For the special case where $f$ is linear, Theorem~\ref{struc_sol_convex} expresses several recent results in \cite[Theorems~3.2--3.4]{Luan_2016} and \cite[Theorem~4.5]{Luan_Yen_2015}.
\end{remark}

The following example is designed as an illustration for Theorem~\ref{struc_sol_convex}.

\begin{example}\label{Ex1}
Let $X=C[-1, 1]$ be the Banach space of continuous real-valued functions defined on $[-1, 1]$ with the norm given by $||x||=\max\limits_{t \in [-1, 1]} |x(t)|.$ The Riesz representation theorem (see, e.g., \cite[Theorem~6,~p.~374]{Kolmogorov_Fomin_1975} and \cite[Theorem~1, p.~113]{Luenberger_1969}) asserts that the dual space of $X$ is $X^*=NBV[-1, 1]$, the \textit{normalized space of functions of bounded variation} on $[-1, 1]$, i.e., functions $y: [-1, 1] \rightarrow \mathbb{R}$ of bounded variation, $y(-1)=0$, and $y(\cdot)$ is continuous from the right at every point of $(-1, 1)$. We define two elements $x^*_1,\, x^*_2 \in X^*$ by setting $$	\langle x^*_i, x \rangle =\int\limits_{-1}^1 t^i x(t)dt\quad (i=1,2),$$ where the integrals are Riemannian. They equal respectively to the Riemann-Stieltjes integrals (see \cite[p.~367]{Kolmogorov_Fomin_1975})
$\int\limits_{-1}^1 x(t)dy_i(t)$, $i=1,2$, which are given by the $C^1$-smooth functions $y_i(t)=\int\limits_{-1}^t \tau^i d \tau$ for $i=1,2$. Now, we consider~$X$ with the \textit{weak topology}. Then $X$ is a locally convex Hausdorff topological vector space whose topology is much weaker than the norm topology. It is clear that $X_0:={\rm ker}\,x^*_1 \cap {\rm ker}\,x^*_2$ is a closed linear subspace of~$X$. Let 
\begin{equation*}
e_1(t)=\begin{cases}
0 & \text{if} \ \, t \in [-1,0]\\
-60t^2+48t   & \text{if} \ \, t \in [0,1], \\ 
\end{cases}
\end{equation*}
 and
\begin{equation*}
e_2(t)=\begin{cases}
0 & \text{if} \ \, t \in [-1,0]\\
80t^2-60t   & \text{if} \ \, t \in [0,1]. \\ 
\end{cases}
\end{equation*}
We have $e_1, e_2\in X$, $\langle x^*_1,e_1 \rangle= \langle x^*_2,e_2 \rangle=1,$ and $\langle x^*_1,e_2 \rangle= \langle x^*_2,e_1 \rangle=0$. For any $x \in X$, put $t_i=\langle x^*_i,x \rangle$ for $i=1,2$, and observe that the vector $x_0:=x-t_1e_1-t_2e_2$ belongs to $X_0$. Conversely, if $x=x_0+t_1e_1+t_2e_2$, with $x_0 \in X_0$ and $t_1, t_2 \in \mathbb{R}$, then
$$\langle x^*_i, x \rangle = \langle x^*_i, x_0 \rangle + t_1\langle x^*_i, e_1 \rangle+ t_2\langle x^*_i, e_2 \rangle=t_i \quad(i=1,2).$$ Therefore, for any $x \in X$, there exists a unique element $(x_0, t_1, t_2) \in X_0 \times \mathbb{R} \times \mathbb{R}$ such that $x=x_0+t_1e_1+t_2e_2$.  

Let $Y=X$, $y^*_i=x^*_i$ for $i=1,2,$ and $K:=\left\{y \in Y \mid \langle y^*_i, y \rangle \leq 0, \, i=1,2\right\}.$ Note that $K$ is a polyhedral convex cone. Clearly,
\begin{equation}\label{ex_rep_K}
K=\left\{x_0 + t_1 e_1 + t_2 e_2 \mid x_0 \in X_0,\; t_i \leq 0,\; i=1,2\right\}.
\end{equation}
By Lemma~\ref{reps_intK}, ${\rm int}\,K=\left\{y \in Y \mid \langle y^*_i, y \rangle < 0, \, i=1,2\right\}$, and  
\begin{equation*}
\begin{aligned}
K\setminus{\ell(K)}&=\left\{y \in Y \mid \langle y^*_1, y \rangle \leq 0, \, \langle y^*_2, y \rangle < 0\right\} \\
&\hspace*{2cm} \cup \left\{y \in Y \mid \langle y^*_1, y \rangle < 0, \, \langle y^*_2, y \rangle \leq 0\right\}.
\end{aligned}
\end{equation*}
An easy computation shows that
\begin{equation}\label{ex_rep_intK}
{\rm int}\,K=\left\{x_0 + t_1 e_1 + t_2 e_2 \mid x_0 \in X_0,\,  t_1 < 0, \, t_2 < 0\right\}
\end{equation}
and 
\begin{equation}\label{ex_rep_diffK}
K\setminus{\ell(K)}=\left\{x_0 + t_1 e_1 + t_2 e_2 \mid x_0 \in X_0,\, t_1 \leq 0, \, t_2 \leq 0, \, t_1+t_2 <0\right\}.
\end{equation}
Given any $e_0 \in X_0$ and put $L=\left\{x \in X \mid x(t)=e_0(t),\, t \in [-1,0] \right\}$. Clearly, $L$ is a closed affine subspace of $X$. Therefore, the set $$D:=\left\{x \in L \mid \langle x^*_1, x \rangle \leq 0, \; \langle x^*_2, x \rangle \leq 1\right\}$$ is generalized polyhedral convex. Observe that $e_0 + e_2 \in D$ and $D$ is not a polyhedral convex set in $X$. 

Let $X_2$ be the linear subspace of $X$ generated by $e_2$. It is clear that ${\rm dim}\,X_2=1$, $X={\rm ker}\,x^*_2 + X_2$, and ${\rm ker}\,x^*_2 \cap X_2 =\{0\}$. Let  $$\pi: X \rightarrow X/{\rm ker}\,x^*_2, \quad \pi(x)= x+{\rm ker}\,x^*_2 \quad (x \in X),$$
be the canonical projection from $X$ on the quotient space $X/{\rm ker}\,x^*_2$. According to \cite[Theorem~1.41(a)]{Rudin_1991}, the linear mapping $\pi$ is continuous. Since the operator 
$$\Phi: X/{\rm ker}\,x^*_2 \rightarrow X_2, \quad x+{\rm ker}\,x^*_2 \mapsto x \quad (x \in X_2),$$ is a linear bijective mapping, $\Phi$ is a homeomorphism by Lemma~\ref{cont_proj_LY_2015}. So, $$\Phi \circ \pi: X \rightarrow X_2$$ is a linear continuous mapping. Define the map $\varphi: X_2 \rightarrow Y$ by $\varphi(te_2)=te_1$ for all $t \in \mathbb{R}$. In accordance with Theorem~3.4 of \cite[p.~22]{Schaefer_1971}, since ${\rm dim}\,X_2=1$, the linear mapping $\varphi$ is continuous. Set $T=\varphi \circ \Phi \circ \pi$, and observe that $T: X \rightarrow X$  is a linear continuous mapping. It is easy to check that if $x=x_0+t_1e_1+t_2e_2$ with $x_0 \in X_0$ and $t_1, t_2 \in \mathbb{R}$, then $Tx=t_2e_1$. Moreover, ${\rm ker}\,T={\rm ker}\,x^*_2$.

Obviously, the space $X$ is the union of the pcs $P_1:=\{x \in X \mid \langle x^*_2, x \rangle \geq 0\}$ and $P_2:=\{x \in X \mid \langle x^*_2, x \rangle \leq 0 \}.$ We see at once that if $x=x_0+t_1e_1+t_2e_2$, where $x_0 \in X_0$ and $t_1, t_2 \in \mathbb{R}$, then $x$ belongs to $P_1$ (resp., belongs to $P_2$) if and only if $t_2 \geq 0$ (resp., $t_2 \leq 0$). Let $f: X \rightarrow Y$ be given by
	\begin{equation*}
	f(x)=\begin{cases}
	x-Tx & {\rm if} \; x \in P_1\\
	x+Tx & {\rm if} \;  x \in P_2.
	\end{cases}
	\end{equation*}
Since $x-Tx=X+Tx=x$ for any $x \in P_1 \cap P_2={\rm ker}\,T$, the values of $f$ is well defined. Moreover, $f$ is a piecewise linear vector-valued function.

\medskip
\noindent
{\sc Claim 1.} \textit{The mapping $f$ is a $K$-function on $X$.}

\smallskip	
Indeed, take any vectors $x_1, \, x_2 \in X$, a number $\lambda \in [0,1]$, and set $$y:=(1-\lambda)f(x_1) + \lambda f(x_2) - f ((1-\lambda) x_1 + \lambda x_2).$$ If there exists an index $i \in \{1,2\}$ such that $x_1, x_2 \in P_i$, then $y=0$ because $f$ is linear on~$P_i$; so $y \in K$. If $x_1 \in P_1$ and $x_2 \in P_2$, then $x_i=x_{i,0}+t_{i,1}e_1+t_{i,2}e_2$, $i=1,2$, with $x_{i,0} \in X_0$, $t_{i,1} \in \mathbb{R}$, $i=1,2$, and $t_{1,2} \geq 0 \geq t_{2,2}$. Therefore, $Tx_i=t_{i,2} e_1$ for $i=1,2$. In the case where $(1-\lambda) x_1 + \lambda x_2 \in P_1$, we obtain
	\begin{equation*}
	\begin{aligned}
	y&=(1-\lambda)(x_1-Tx_1) + \lambda (x_2+Tx_2)\\
	&\hspace*{3cm} -  ((1-\lambda) x_1 + \lambda x_2)+T((1-\lambda) x_1 + \lambda x_2)\\
	&=-(1-\lambda)Tx_1 + \lambda Tx_2+(1-\lambda)Tx_1 + \lambda Tx_2\\
	&=2\lambda (Tx_2)=2\lambda t_{2,2} e_1.
	\end{aligned}
	\end{equation*}
Since $2\lambda t_{2,2} \leq 0$, one has $y \in K$ by \eqref{ex_rep_K}. If $(1-\lambda) x_1 + \lambda x_2 \in P_2$, then
\begin{equation*}
\begin{aligned}
y&=(1-\lambda)(x_1-Tx_1) + \lambda (x_2+Tx_2)\\
&\hspace*{3cm} -  ((1-\lambda) x_1 + \lambda x_2)-T((1-\lambda) x_1 + \lambda x_2)\\
&=-(1-\lambda)Tx_1 + \lambda Tx_2-(1-\lambda)Tx_1 - \lambda Tx_2\\
&=-2(1-\lambda)Tx_1=-2(1-\lambda) t_{1,2} e_1.
\end{aligned}
\end{equation*}
From \eqref{ex_rep_K}, since $-2(1-\lambda) t_{1,2} \leq 0$, one gets $y \in K$. With $x_1 \in P_2$ and $x_2 \in P_1$, a similar conclusion is obtained. It follows that $f$ is a $K$-function on $X$. 

\medskip
\noindent
{\sc Claim 2.} \textit{It holds that}
\begin{equation}\label{ex1_sol_convex}
	{\rm Sol(VP)}=\left\{u \in L \mid \langle x^*_1, u \rangle = 0, \, 0 \leq \langle x^*_2, u \rangle \leq 1\right\}.
\end{equation}

\smallskip
First, take any $u \in {\rm Sol(VP)}$. Let $u_0 \in X_0$ and $t_1, t_2 \in \mathbb{R}$ be such that $u=u_0+t_1e_1+t_2e_2$. Since $u \in D$, $t_1 \leq 0$ and $t_2 \leq 1$. Moreover, $u_0(t)=e_0(t)$ for all $t \in [-1,0]$; so, $u_0 \in D$. If $t_2 <0$, then $u \in P_2$. Therefore, $$f(u)=u+Tu=u_0+(t_1+t_2)e_1+t_2e_2.$$ Since $f(u)-f(u_0)=(t_1+t_2)e_1+t_2e_2$ with $t_1+t_2 < 0$ and $t_2<0$, by~\eqref{ex_rep_diffK}, $f(u)-f(u_0) \in K\setminus \ell(K)$. The fact contradicts the assumption $u \in {\rm Sol(VP)}$. This clearly forces $0 \leq t_2 \leq 1$; hence $u \in P_1$. If $t_1 <0$, then we choose $x=u_0+t_2e_2$. Since $x \in D \cap P_1$,
	\begin{equation*}
	\begin{aligned}
	f(u)-f(x)&=(u-Tu)-(x-Tx)\\
	&=(u_0+t_1e_1+t_2e_2-t_2e_1)-(u_0+t_2e_2-t_2e_1)=t_1e_1.\\
	\end{aligned}
	\end{equation*}
Since $t_1<0$, by \eqref{ex_rep_diffK}, $f(u)-f(x) \in K\setminus \ell(K)$. This inclusion contradicts the assumption $u \in {\rm Sol(VP)}$. We thus get $t_1=0$ and $0 \leq t_2 \leq 1$. Consequently, $u \in S$, where $S$ is the set on the right-hand side of \eqref{ex1_sol_convex}. We have proved that ${\rm Sol(VP)} \subset S$. To obtain the opposite inclusion, take any $u \in S$. Let us show that $f(u)-f(x) \notin K\setminus \ell(K)$ for all $x \in D$. Suppose that $u=u_0+t_1e_1+t_2e_2$ with $u_0 \in X_0$ and $t_1, t_2 \in \mathbb{R}$. Of course, $t_1=0$ and $0 \leq t_2 \leq 1$. Since $u \in P_1$, $$f(u)=u-Tu=u_0-t_2e_1+t_2e_2.$$ 
For any $x \in D$, there exist a vector $x_0 \in X_0$, numbers $\tau_1 \leq 0$ and $\tau_2 \leq 1$ satisfying $x=x_0+\tau_1e_1+\tau_2e_2$. 

If $\tau_2 <0$, then $x \in P_2$; so $f(x)=x+Tx=x_0+(\tau_1+\tau_2)e_1+\tau_2e_2$. According to \eqref{ex_rep_diffK}, since $f(u)-f(x)=(u_0-x_0)+ (-t_2-\tau_1-\tau_2)e_1+(t_2 -\tau_2) e_2$ with $t_2-\tau_2 >0$, we can assert that $f(u)-f(x)  \notin K\setminus \ell(K)$.

If $\tau_2 \geq 0$ then $x \in P_1$. Since $f(x)=x-Tx=x_0+(\tau_1-\tau_2)e_1+\tau_2e_2$, 
\begin{equation*}
f(u)-f(x)=(u_0-x_0)+ (-t_2-\tau_1+\tau_2)e_1+(t_2-\tau_2) e_2.
\end{equation*}
As $(-t_2-\tau_1+\tau_2) + (t_2-\tau_2)=-\tau_1 \geq 0,$ one gets $f(u)-f(x) \notin K\setminus \ell(K)$ by~\eqref{ex_rep_diffK}. From what that has already been said, we obtain $u \in {\rm Sol}{\rm (VP)}$. We have proved that ${\rm Sol(VP)}=S$.

\medskip
\noindent
{\sc Claim 3.} \textit{It holds that}
\begin{equation}\label{ex1_wsol_convex}
	\begin{aligned}
	{\rm Sol}^w{\rm (VP)}&=\left\{u \in L \mid \langle x^*_1, u \rangle = 0, \, 0 \leq \langle x^*_2, u \rangle \leq 1\right\}\\
	&\hspace*{2cm} \cup \left\{u \in L \mid \langle x^*_1, u \rangle \leq 0, \, \langle x^*_2, u \rangle =1  \right\}.
	\end{aligned}
\end{equation}

First, take any $u \in 	{\rm Sol}^w{\rm (VP)}$. Suppose that $u=u_0+t_1e_1+t_2e_2$ with $u_0 \in X_0$ and $t_1, t_2 \in \mathbb{R}$. Since $u \in D$, $t_1 \leq 0$ and $t_2 \leq 1$. Moreover, $u_0(t)=e_0(t)$ for all $t \in [-1,0]$. Therefore, $u_0 \in D$. If $t_2 <0$ then $u \in P_2$; so $$f(u)=u+Tu=u_0+(t_1+t_2)e_1+t_2e_2.$$ Since $f(u)-f(u_0)=(t_1+t_2)e_1+t_2e_2$ with $t_1+t_2 < 0$ and $t_2 < 0$, by \eqref{ex_rep_intK}, $f(u)-f(u_0) \in {\rm int}\,K$. This contradicts the assumption $u \in {\rm Sol}^w{\rm (VP)}$. We thus get $0 \leq t_2 \leq 1$. Therefore, $u \in P_1$ and $$f(u)=u-Tu=u_0+(t_1-t_2)e_1+t_2e_2.$$ If $t_1 <0$ and $t_2 < 1$, one can find a positive number $\varepsilon$ such that $t_1+\varepsilon <0$ and $t_2+\varepsilon <1$. Take $x=u_0+(t_2+\varepsilon)e_2$, and observe that $x \in D \cap P_1$. Since $f(x)=x-Tx=u_0-(t_2+\varepsilon)e_1+(t_2+\varepsilon)e_2$,
\begin{equation*}
f(u)-f(x)=(t_1+\varepsilon)e_1+(-\varepsilon)e_2.
\end{equation*}
As $t_1+\varepsilon <0$ and $-\varepsilon<0$, formula \eqref{ex_rep_intK} shows that $f(u)-f(x) \in {\rm int}\,K$, which is impossible because $u \in {\rm Sol}^w{\rm (VP)}$. We thus get $t_1=0$ or $t_2=1$. Consequently, $u$ belongs to $S^w$, where $S^w$ is the set on the right-hand side of~\eqref{ex1_wsol_convex}. We have proved that ${\rm Sol}^w{\rm (VP)} \subset S^w$. To obtain the opposite inclusion, take any $u \in S^w$. Let $u_0 \in X_0$ and $t_1, t_2 \in \mathbb{R}$ be such that $u=u_0+t_1e_1+t_2e_2$. 

\smallskip\noindent
{\it Case 1:} $t_1=0$ and $0 \leq t_2 \leq 1$. Clearly, $\langle x^*_1, u \rangle = 0$ and $0 \leq \langle x^*_2, u \rangle \leq 1$. According to \eqref{ex1_sol_convex}, $u \in {\rm Sol}{\rm (VP)}$; hence $u \in {\rm Sol}^w{\rm (VP)}$.

\smallskip\noindent
{\it Case 2:} $t_1<0$ and $t_2=1$. Since $u \in P_1$, one has $$f(u)=u-Tu=u_0+(t_1-1)e_1+e_2.$$ For any $x \in D$, there are $x_0 \in X_0$ and $\tau_1, \tau_2 \in \mathbb{R}$ satisfying $x=x_0+\tau_1e_1+\tau_2e_2$. Clearly, $\tau_1 \leq 0$ and $\tau_2 \leq 1$. If $\tau_2 <0$, then $x \in P_2$ and $$f(x)=x+Tx=x_0+(\tau_1+\tau_2)e_1+\tau_2e_2.$$ According to \eqref{ex_rep_intK}, since $f(u)-f(x)=(u_0-x_0)+ (t_1-1-\tau_1-\tau_2)e_1+(1 -\tau_2) e_2$ with $1-\tau_2 >0$, $f(u)-f(x) \notin {\rm int}\,K$. If $0 \leq \tau_2 \leq 1$, then $x \in P_1$ and $$f(x)=x-Tx=x_0+(\tau_1-\tau_2)e_1+\tau_2e_2.$$ Since $f(u)-f(x)=(u_0-x_0)+ (t_1-1-\tau_1+\tau_2)e_1+(1 -\tau_2) e_2$ with $1-\tau_2 \geq 0$, formula \eqref{ex_rep_intK} gives $f(u)-f(x)  \notin {\rm int}\,K$. It follows that $f(u)-f(x) \notin {\rm int}\,K$ for all $x \in D$. Hence, $u \in {\rm Sol}^w{\rm (VP)}$. We have proved that ${\rm Sol}^w{\rm (VP)}=S^w$.

\medskip
Observe that ${\rm Sol}^w{\rm (VP)} \neq {\rm Sol}{\rm (VP)}$. Indeed, let $e_3$ in $X$ be given by 
\begin{equation*}
e_3(t)=\begin{cases}
0 & \text{if} \ \, t \in [-1,0]\\
155t^2-120t   & \text{if} \ \, t \in [0,1]. 
\end{cases}
\end{equation*}
Since $e_0+e_3 \in L$, $\langle x^*_1, e_0+e_3 \rangle = - \frac{5}{4}$, and $\langle x^*_2, e_0+e_3 \rangle = 1$, we can see that $e_0+e_3$ belongs to ${\rm Sol}^w{\rm (VP)}$ but $e_0+e_3 \notin {\rm Sol}{\rm (VP)}$.
\end{example}
 
\section{Structure of the Solution Sets in the Nonconvex Case}
In Theorem~\ref{struc_sol_convex}, the assumption $f$ is a $K$-function on $D$ cannot be dropped (see \cite[Example~2.1]{Yang_Yen_2010} for an example about the efficient solution set, \cite[p.~1252]{ZhengYang_2008} for an example about the weakly efficient solution set). For the case where $X,\, Y$ are normed spaces, $Y$ is of finite dimension, $K \subset Y$ is a pointed cone, and $D\subset X$ is a polyhedral convex set, the efficient solution set of {\rm (VP)} is shown to be the union of finitely many semi-closed polyhedral convex sets (see \cite[Theorem~2.1]{Yang_Yen_2010}). This result can be extended to the lcHtvs setting we are studying as follows.

\begin{theorem}\label{struc_sol_nonconvex} The efficient solution set of {\rm (VP)} is the union of finitely many semi-closed generalized polyhedral convex sets.
\end{theorem}

\noindent{\it Proof} \, We may assume, without loss of generality, that $D$ is a nonempty set. 

\medskip
\noindent
{\sc Claim 1.} 	\textit{There exists a closed linear subspace $X_1$ of $X$ such that $D$ is a polyhedral convex set in $X_1$.}

\smallskip
Suppose that $D$ is given by \eqref{eq_def_gpcs_2}. Set $X_0:={\rm ker}\,A$ and observe that $X_0$ is a closed linear subspace of $X$. 

If $z=0$, then $D=\{x \in X_0 \mid \langle x_i^*, x \rangle \leq \alpha_i,\, i=1,\dots,p\}$. For $i=1,\dots,p$, let $x^*_{i,0}$ be the restriction to $X_0$ of the functional $x^*_i$. Clearly, $x^*_{i,0} \in X^*_0$ for $i=1,\dots,p$. Since $$D=\{x \in X_0 \mid \langle x_{i,0}^*, x \rangle \leq \alpha_i,\, i=1,\dots,p\},$$ we can assert that $D$ is a polyhedral convex set in $X_0$.

If $z \neq 0$, then fix any vector $\bar x \in D$. Clearly, $\bar x$ does not belong to $X_0$.  Set $X_1:=X_0 + \{t \bar x \mid t \in \mathbb{R}\}.$ Since $\{t \bar x \mid t \in \mathbb{R}\}$ is a linear subspace of one-dimension, the linear subspace $X_1$ is closed by \cite[Theorem~1.42]{Rudin_1991}. According to \cite[Theorem~3.5]{Rudin_1991}, there exists $x_{1,0}^* \in X_1^*$ such that $\langle x_{1,0}^*,  \bar x \rangle =1$ and $\langle x_{1,0}^*,  x \rangle =0$ for all $x \in X_0$. It is not difficult to show that $X_0={\rm ker}\,x_{1,0}^*$. For $i=1,\dots,p$, let $x^*_{1,i}$ be the restriction to $X_1$ of the functional $x^*_i$. Of course, $x^*_{1,i} \in X^*_1$ for $i=1,\dots,p$. Set
$$D':=\{x \in X_1 \mid \langle x_{1,0}^*, x \rangle =1, \langle x_{1,i}^*, x \rangle \leq \alpha_i, i=1,\dots,p\},$$ and observe that $D'$ is a polyhedral convex set in $X_1$. Let us show that $D=D'$. To obtain the inclusion $D \subset D'$, take any $x \in D$. As $x_0:=x-\bar x$ belongs to~$X_0$, one has $x \in X_1$. Therefore, $$\langle x_{1,0}^*, x \rangle= \langle x_{1,0}^*, \bar x \rangle + \langle x_{1,0}^*, x_0 \rangle=\langle x_{1,0}^*, \bar x \rangle =1.$$ For each $i=1,\dots,p$, since $\langle x^*_{1,i}, x \rangle = \langle x^*_{i}, x \rangle$, one gets $\langle x_{1,i}^*, x \rangle \leq \alpha_i$. It follows that $x \in D'$.  We have proved that $D \subset D'$. To obtain the opposite inclusion, take any $x \in D'$. Let $x_0 \in X_0$ and $t \in \mathbb{R}$ be such that $x=x_0+t\bar x$. Then $$\langle x_{1,0}^*, x \rangle =\langle x_{1,0}^*, x_0 \rangle + t \langle x_{1,0}^*, \bar x \rangle=t.$$ Since $x \in D'$, $\langle x_{1,0}^*, x \rangle =1$; so $t=1$. Therefore, $$Ax=A(x_0+\bar x)=Ax_0+A\bar x =z.$$ For each $i=1,\dots,p$, the inequality $\langle x_{1,i}^*, x \rangle \leq \alpha_i$ implies $\langle x^*_{i}, x \rangle  \leq \alpha_i$. The inclusion $D' \subset D$ has been proved. Thus $D=D'$; hence $D$ is a polyhedral convex set in $X_1$.  

\medskip
By Claim 1, we may assume that $D$ is a pcs in $X_1$, where $X_1$ is a closed linear subspace of $X$. For each $k=1,\dots,m,$ set $P_{1,k}:=P_k \cap X_1$ and observe that $P_{1,k}$ is a pcs in $X_1$. Of course, $X_1=\bigcup\limits_{k=1}^m P_{1,k}$. Let $$\pi_1: Y \rightarrow Y/Y_0, \quad \pi_1(y)= y+Y_0 \quad (y \in Y),$$ be the canonical projection from $Y$ on the quotient space $Y/Y_0$.  By \cite[Theorem~1.41(a)]{Rudin_1991}, $\pi_1$ is a linear continuous mapping. Since the operator $$\Phi_1: Y/Y_0 \rightarrow Y_1, \quad y +Y_0 \mapsto y \quad (y \in Y_1),$$ is a linear bijective mapping, $\Phi_1$ is a homeomorphism by Lemma~\ref{cont_proj_LY_2015}. So, the operator $\pi:=\Phi_1\circ \pi_1: Y \rightarrow Y_1$ is linear and continuous. 

\medskip\noindent
{\sc Claim 2.} \textit{For any $y \in Y$, we have $\pi (y) \in K_1 \setminus \{0\}$ if and only if $y \in K \setminus \ell(K)$.}

\smallskip
Indeed, suppose that $y \in K \setminus \ell(K)$. By Lemma~\ref{decomp_K_diff_lK}, one can find $y_0 \in Y_0$ and $y_1 \in K_1 \setminus \{0\}$ such that $y=y_0 + y_1$. Then $\pi(y)=y_1 \in K_1 \setminus \{0\}$. Now, suppose that $\pi(y)  \in K_1 \setminus \{0\}$, i.e., there exists $y_1 \in  K_1 \setminus \{0\}$ satisfying $\pi(y)=y_1$. This implies that $y-y_1 \in Y_0$; hence $y \in y_1 + Y_0 \subset  K_1 \setminus \{0\}+Y_0$. In accordance with Lemma~\ref{decomp_K_diff_lK}, $y \in K \setminus \ell(K)$. 

\medskip
Let $f_1$ be the restriction to $X_1$ of the mapping $\pi \circ f$. Clearly, $f_1$ is a piecewise linear vector-valued function from $X_1$ to $Y_1$. Let us consider a piecewise linear vector optimization problem   
\begin{equation*}
{\rm (VP_1)} \qquad \quad {\rm Min}_{K_1} \big\{ f_1(x) \mid x \in D\big\}.
\end{equation*}

\smallskip\noindent
{\sc Claim 3.} \textit{It holds that ${\rm Sol(VP)}={\rm Sol(VP_1)}$.}

\smallskip
First, to show that ${\rm Sol(VP)} \subset {\rm Sol(VP_1)}$, we suppose the contrary: There exists $u \in {\rm Sol(VP)}$ not belonging to ${\rm Sol(VP_1)}$. Then one can find $x \in D$ such that $f_1(u)-f_1(x) \in K_1 \setminus \{0\}.$ Since $\pi (f(u)-f(x)) \in K_1 \setminus \{0\}$, by Claim~2, one has $f(u)-f(x) \in K \setminus \ell(K)$. This contradicts the assumption $u \in {\rm Sol(VP)}$. 

Now, to obtain the inclusion ${\rm Sol(VP_1)} \subset {\rm Sol(VP)}$, take any $u \notin {\rm Sol(VP)}$. Then there exists $x \in D$ such that $f(u)-f(x) \in K \setminus \ell(K)$. Combining this with Claim 2, one gets $\pi (f(u)-f(x)) \in K_1 \setminus \{0\}$, i.e., $f_1(u)-f_1(x) \in K_1 \setminus \{0\}.$ Therefore, $u \notin {\rm Sol(VP_1)} $.  We have thus proved that ${\rm Sol(VP)}={\rm Sol(VP_1)}$. 

\medskip
Since $Y_1$ is finite-dimensional, $K_1$ is a pointed cone, $D$ is a polyhedral convex set in $X_1$, arguing similarly as in the proof of \cite[Theorem~2.1]{Yang_Yen_2010}, we can assert that ${\rm Sol(VP_1)}$ is the union of finitely many semi-closed polyhedral convex sets in $X_1$. As ${\rm Sol(VP)}={\rm Sol(VP_1)}$ by Claim 3, the assertion of the theorem has been proved.  $\hfill\Box$ 

\medskip
The next result is a generalization of \cite[Theorem~3.1]{ZhengYang_2008}.
\begin{theorem}\label{struc_wsol_nonconvex} If ${\rm int}\,K$ is nonempty, then the weakly efficient solution set of~{\rm (VP)} is the union of finitely many generalized polyhedral convex sets.
\end{theorem}
\noindent{\it Proof} \, We may assume, without loss of generality, that $D$ is nonempty. Set 
$$H_{j}:=\{y \in Y \mid \langle y^*_{j}, y \rangle \leq 0 \} \quad (j=1,\dots,q).$$ Invoking Lemma~\ref{reps_intK}, we have ${\rm int}\,H_j=\{y \in Y \mid \langle y^*_{j}, y \rangle < 0 \}$ for $j=1,\dots,q,$ and ${\rm int}\,K=\bigcap\limits_{j=1}^q {\rm int}\,H_j$. By Claim~1 in the proof of Theorem~\ref{struc_sol_convex}, one sees $M_k+K$ is a polyhedral convex set in $Y$, where $M_k:=f(D \cap P_k)$ for all $k=1,\dots,m.$ Then, one can find $y^*_{k,1}, \dots, y^*_{k,\ell_k}$ in $Y^*$, $\beta_{k,1}, \dots, \beta_{k,\ell_k}$ in~$\mathbb{R}$ such that  $$M_k+K=\bigcap\limits_{i=1}^{\ell_k}H_{k,i},$$ with $H_{k,i}:=\{y \in Y \mid \langle y^*_{k,i}, y \rangle \leq \beta_{k,i} \}.$ Put $Q=f(D)+K$ and observe that $Q=\bigcup\limits_{k=1}^{m}(M_k+K)$. Since $E^w(Q|K) = Q \setminus (Q + {\rm int}\,K),$ one has
\begin{equation}\label{proof_weak_sol_1}
\begin{aligned}
E^w(Q|K) & =Q\setminus \Big( \bigcup\limits_{k=1}^{m}  (M_{k}+K)+ {\rm int}\,K\Big)\\
&=\bigcap\limits_{k=1}^{m} \Big( Q\setminus  (M_{k}+K+ {\rm int}\,K)\Big)\\
&=\bigcap\limits_{k=1}^{m} \Big(Q\setminus \big( \bigcap\limits_{i=1}^{\ell_{k}}H_{k,i}+ {\rm int}\,K\big)\Big)\\
&=\bigcap\limits_{k=1}^{m} \bigcup\limits_{i=1}^{\ell_{k}} \Big( Q\setminus \big(H_{k,i}+{\rm int}\,K) \Big).
\end{aligned}
\end{equation}
For any $k \in \{1,\dots,m\}$ and $i \in \{1,\dots,\ell_{k}\}$, 
\begin{equation}\label{proof_weak_sol_2}
\begin{aligned}
Q\setminus \big(H_{k,i}+{\rm int}\,K\big)
&=Q\setminus \left( \bigcap\limits_{j=1}^q \big(H_{k,i}+{\rm int}\,H_j\big) \right)\\
&=\bigcup\limits_{j=1}^q \Big( Q \setminus \big(H_{k,i}+{\rm int}\,H_j\big) \Big)\\
&=\bigcup\limits_{j=1}^q \bigcup\limits_{k_1=1}^m \Big( (M_{k_1}+K) \setminus \big(H_{k,i}+{\rm int}\,H_j\big)  \Big). 
\end{aligned}
\end{equation}

\medskip
Observe that, \textit{for any $k,\, k_1 \in \{1,\dots,m\}$, $i \in \{1,\dots,\ell_{k}\}$, and $j \in \{1,\dots,q\}$, the set $(M_{k_1}+K) \setminus \big(H_{k,i}+{\rm int}\,H_j\big)$ is polyhedral convex.}

\smallskip
Indeed, let us show that there exist $y^*_{k,i,j} \in Y^*$ and $\beta_{k,i,j} \in \mathbb{R}$ such that
\begin{equation}\label{proof_weak_sol_3}
H_{k,i}+{\rm int}\,H_j=\{y \in Y \mid \langle y^*_{k, i, j}, y \rangle < \beta_{k,i,j} \}.
\end{equation}
First, consider the case where $y^*_{k,i}=0$. If $\beta_{k,i} < 0$, then $H_{k,i}$ is empty. So, one can choose $y^*_{k, i, j}=0$ and $\beta_{k,i,j}=\beta_{k,i}$ because $H_{k,i}+{\rm int}\,H_j=\emptyset$.  If $\beta_{k,i} \geq 0$, then $H_{k,i}=Y$. Since $H_{k,i}+{\rm int}\,H_j=Y$, \eqref{proof_weak_sol_3} is fulfilled with $y^*_{k, i, j}=0$ and $\beta_{k,i,j}=\beta_{k,i}.$ Now, we consider the case where $y^*_{k,i} \neq 0$. One can find $\bar y_{k,i}, w_{k,i}$ in $Y$ satisfying $\langle  y^*_{k,i}, \bar y_{k,i} \rangle=\beta_{k,i}$ and $\langle  y^*_{k,i}, w_{k,i} \rangle=1$. It is not difficult to show that $$H_{k,i}=\bar y_{k,i} + \{t_{k,i} w_{k,i} \mid t_{k,i}\leq 0\}+{\rm ker}\,y^*_{k,i}.$$
Since $y^*_j \neq 0$, there exists $w_{j} \in Y$ such that $\langle  y^*_{j}, w_{j} \rangle=1$. We see at once that $H_{j}=\{t_j w_{j} \mid t_j \leq 0\}+{\rm ker}\,y^*_{j}$ and ${\rm int}\,H_{j}=\{t_j w_{j} \mid t_j < 0\}+{\rm ker}\,y^*_{j}$. It follows that
\begin{equation}\label{rep_sum_hca}
\begin{aligned}
H_{k,i}+{\rm int}\,H_{j}&=\bar y_{k,i} + \{t_{k,i} w_{k,i} \mid t_{k,i}\leq 0\}+\{t_j w_{j} \mid t_j < 0\}\\
&\hspace*{5cm}+{\rm ker}\,y^*_{k,i}+{\rm ker}\,y^*_{j}.
\end{aligned}
\end{equation}

If ${\rm ker}\,y^*_{k,i}\neq {\rm ker}\,y^*_{j}$, then ${\rm ker}\,y^*_{k,i}+{\rm ker}\,y^*_{j}=Y$ by ${\rm codim}({\rm ker}\,y^*_{k,i})=1.$ Therefore, \eqref{rep_sum_hca} shows that $H_{k,i}+{\rm int}\,H_j=Y$. So, one can choose $y^*_{k, i, j}=0$ and $\beta_{k,i,j}=1$.    

If ${\rm ker}\,y^*_{k,i}={\rm ker}\,y^*_{j}$, then we take $\lambda_{k,i,j}=\langle y^*_j, w_{k,i} \rangle$. For every $y \in Y$, put $t_{k,i}=\langle y^*_{k,i}, y\rangle$. Clearly, the vector $y_0:=y-t_{k,i}w_{k,i}$ belongs to ${\rm ker}\,y^*_{k,i}$. Therefore, 
\begin{equation*}
\begin{aligned}
\langle y^*_j-\lambda_{k,i,j}y^*_{k, i},y \rangle &= \langle y^*_j-\lambda_{k,i,j}y^*_{k, i},y_0 \rangle+\langle y^*_j-\lambda_{k,i,j}y^*_{k, i},t_{k,i}w_{k,i} \rangle\\
&= \langle y^*_j,y_0 \rangle - \lambda_{k,i,j} \langle y^*_{k, i},y_0 \rangle+ t_{k,i}\langle  y^*_j-\lambda_{k,i,j}y^*_{k, i},w_{k,i} \rangle\\
&= t_{k,i} \left(\langle  y^*_j,w_{k,i} \rangle - \lambda_{k,i,j}\langle  y^*_{k, i},w_{k,i} \rangle \right)\\
&=t_{k,i} \big(\lambda_{k,i,j} - \lambda_{k,i,j}\big)=0.\\
\end{aligned}
\end{equation*} 
We thus get $y^*_j=\lambda_{k,i,j}y^*_{k, i, j}$. Since $y^*_j \neq 0$, one has $\lambda_{k,i,j} \neq 0$. If $\lambda_{k,i,j} > 0$, then ${\rm int}\,H_{j}=\{y \in Y \mid \langle y^*_{k, i}, y \rangle < 0\}$. So, $$H_{k,i}+{\rm int}\,H_{j}=\{y \in Y \mid \langle y^*_{k, i}, y \rangle < \beta_{k,i}\}.$$ Of course, the formula \eqref{proof_weak_sol_3} is fulfilled with $y^*_{k, i, j}=y^*_{k, i}$ and $\beta_{k,i,j}=\beta_{k,i}.$  If $\lambda_{k,i,j} < 0$, then ${\rm int}\,H_{j}=\{y \in Y \mid \langle y^*_{k, i}, y \rangle > 0\}$. Therefore, $H_{k,i}+{\rm int}\,H_{j}=Y$. Hence, one can choose $y^*_{k, i, j}=0$ and $\beta_{k,i,j}=1$.  

From \eqref{proof_weak_sol_3} we see that $$(M_{k_1}+K) \setminus \big(H_{k,i}+{\rm int}\,H_j\big)=(M_{k_1}+K) \cap \{y \in Y \mid \langle y^*_{k, i, j}, y \rangle \geq \beta_{k,i,j} \}.$$ Therefore, $(M_{k_1}+K) \setminus \big(H_{k,i}+{\rm int}\,H_j\big)$ is a polyhedral convex set in $Y$. 

\medskip
From \eqref{proof_weak_sol_2} it follows that, for all $k \in \{1,\dots,m\}$ and $i \in \{1,\dots,\ell_{k}\}$, $Q\setminus \big(H_{k,i}+{\rm int}\,K\big)$ is the union of finitely many gpcs. Therefore, by \eqref{proof_weak_sol_1}, $E^w(Q|K)$ is the union of finitely many gpcs. Hence, using the same argument for getting Claim 4 in the proof of Theorem~\ref{struc_sol_convex}, we can prove that ${\rm Sol}^w{\rm (VP)}$ is the union of finitely many gpcs.  $\hfill\Box$ 

\begin{remark} For the case where $X,\, Y$ are normed spaces and $D$ is a polyhedral convex set, the result in Theorem~\ref{struc_wsol_nonconvex} is due to Zheng and Yang \cite[Theorem~3.1]{ZhengYang_2008}.
\end{remark}

Let us  consider an illustrative example for Theorems \ref{struc_sol_nonconvex} and \ref{struc_wsol_nonconvex}. 
\begin{example}
Keeping the notations of Example \ref{Ex1}, we redefine the piecewise linear function $f$ by 
\begin{equation*}
f(x)=\begin{cases}
x+Tx & {\rm if} \; x \in P_1\\
x-Tx & {\rm if} \;  x \in P_2.
\end{cases}
\end{equation*}

\medskip\noindent
{\sc Claim 1.} \textit{It holds that}
\begin{equation}\label{ex_sol_nonconvex}
\begin{aligned}
{\rm Sol(VP)}&=\left\{u \in L \mid \langle x^*_1, u \rangle = 0, \, \langle x^*_2, u \rangle = 1\right\}\\
&\hspace*{2cm} \cup \left\{u \in L \mid \langle x^*_1, u \rangle = 0, \, \langle x^*_2, u \rangle < -1  \right\}.
\end{aligned}
\end{equation}

\smallskip
First, to show that ${\rm Sol(VP)} \subset S$, where $S$ is the set on the right-hand side of \eqref{ex_sol_nonconvex}, take any $u \in {\rm Sol(VP)}$. Then one can find a vector $u_0 \in X_0$ and numbers $t_1, t_2$ satisfying $u=u_0+t_1e_1+t_2e_2$. Since $u \in D$, we have $t_1 \leq 0$ and $t_2 \leq 1$. Moreover, $u_0(t)=e_0(t)$ for all $t \in [-1,0]$; so, $u_0 \in D$. The condition $u \in {\rm Sol(VP)}$ yields $f(u)-f(x) \notin K\setminus \ell(K)$ for every $x \in D$. 

If $t_2 \geq 0$, then $u \in P_1$; so $f(u)=u+Tu=u_0+(t_1+t_2)e_1+t_2e_2$. Observe that $x:=u_0+e_2$ belongs to $D \cap P_1$. It is clear that $f(x)=x+Tx=u_0+e_1+e_2$; hence  
\begin{equation*}
f(u)-f(x)=(t_1+t_2-1)e_1+(t_2-1)e_2.
\end{equation*}
Since $t_1+t_2 -1 \leq 0$ and $t_2-1 \leq 0$, \eqref{ex_rep_K}  shows that $f(u)-f(x) \in K$. Combining this with the inclusion $f(u)-f(x) \notin K\setminus \ell(K)$, one gets $f(u)-f(x) \in \ell(K)$. This implies that $t_1+t_2 -1=0$ and $t_2-1=0$, i.e., $t_1=0$ and $t_2=1$. Therefore, $u \in S$.

If $t_2 < 0$, then $u \in P_2$ and $f(u)=u-Tu=u_0+(t_1-t_2)e_1+t_2e_2$. If we take $x=u_0+e_2$, then $x\in D \cap P_1$. As $f(x)=x+Tx=u_0+e_1+e_2$, one has  
\begin{equation*}
f(u)-f(x)=(t_1-t_2-1)e_1+(t_2-1)e_2.
\end{equation*}
Since $f(u)-f(x) \notin K\setminus \ell(K)$, \eqref{ex_rep_K}  shows that $t_1-t_2-1 >0$, by $t_2-1<0$. Therefore, $t_2 < -1$ as $t_1 \leq 0$. If $t_1<0$, then one can find a positive number $\varepsilon$ such that $t_1+\varepsilon <0$ and $t_2+\varepsilon <-1$. Set $x=u_0+(t_2+\varepsilon)e_2$, and observe that $x \in D \cap P_2$. Since $f(x)=x-Tx=u_0-(t_2+\varepsilon)e_1+(t_2+\varepsilon)e_2$,
\begin{equation*}
f(u)-f(x)=(t_1+\varepsilon)e_1+(-\varepsilon)e_2.
\end{equation*}
As $t_1+\varepsilon <0$ and $-\varepsilon<0$, \eqref{ex_rep_diffK} yields $f(u)-f(x) \in K\setminus \ell(K)$. This contradicts the assumption $u \in {\rm Sol}{\rm (VP)}$. We thus get $t_1=0$. Consequently, $u \in S$. 

We have proved that ${\rm Sol(VP)} \subset S$. To obtain the opposite inclusion, take any $u \in S$. Let $u_0 \in X_0$ and $t_1, t_2 \in \mathbb{R}$ be such that $u=u_0+t_1e_1+t_2e_2$. Of course, $t_1=0$. Given any $x \in D$, one can find a vector $x_0 \in X_0$, numbers $\tau_1 \leq 0$ and $\tau_2 \leq 1$ such that $x=x_0+\tau_1e_1+\tau_2e_2$. 

If $t_2=1$, then $u \in P_1$ and $f(u)=u+Tu=u_0+e_1+e_2$. If $0 \leq \tau_2 \leq 1$, then $x \in P_1$ and $f(x)=x_0+(\tau_1+\tau_2)e_1+\tau_2e_2$. Since 
\begin{equation*}
f(u)-f(x)=(u_0-x_0)+(1-\tau_1-\tau_2)e_1+(1-\tau_2)e_2,
\end{equation*}
with $1-\tau_1-\tau_2 \geq 0$ and $1-\tau_2 \geq 0$, \eqref{ex_rep_diffK} shows that $f(u)-f(x) \notin K\setminus \ell(K)$. If $\tau_2 <0$, then $x \in P_2$; so $f(x)=x_0+(\tau_1-\tau_2)e_1+\tau_2e_2$. According to \eqref{ex_rep_diffK}, since $f(u)-f(x)=(u_0-x_0)+(1-\tau_1+\tau_2)e_1+(1-\tau_2)e_2$
with $$(1-\tau_1+\tau_2)+(1-\tau_2)=2-\tau_1 > 0,$$ one gets $f(u)-f(x) \notin K\setminus \ell(K)$. It follows that $f(u)-f(x) \notin K\setminus \ell(K)$ for all $x \in D$; hence, $u \in {\rm Sol}{\rm (VP)}$. 

If $t_2<-1$, then $u \in P_2$. Therefore, $f(u)=u-Tu=u_0-t_2e_1+t_2e_2$. If $0 \leq \tau_2 \leq 1$, then $x \in P_1$ and $f(x)=x_0+(\tau_1+\tau_2)e_1+\tau_2e_2$. As
\begin{equation*}
f(u)-f(x)=(u_0-x_0)+(-t_2-\tau_1-\tau_2)e_1+(t_2-\tau_2)e_2
\end{equation*}
with $-t_2-\tau_1-\tau_2 > 1-\tau_1-\tau_2 \geq 0$, one has $f(u)-f(x) \notin K\setminus \ell(K)$ by \eqref{ex_rep_diffK}. If $\tau_2 <0$, then $x \in P_2$ and $f(x)=x_0+(\tau_1-\tau_2)e_1+\tau_2e_2$. Observe that
\begin{equation*}
f(u)-f(x)=(u_0-x_0)+(-t_2-\tau_1+\tau_2)e_1+(t_2-\tau_2)e_2
\end{equation*}
with $(-t_2-\tau_1+\tau_2)+(t_2-\tau_2)=-\tau_1 \geq 0$. So, $f(u)-f(x) \notin K\setminus \ell(K)$ by~\eqref{ex_rep_diffK}. Therefore, $f(u)-f(x) \notin K\setminus \ell(K)$ for all $x \in D$. Hence, $u \in {\rm Sol}{\rm (VP)}$. We have proved that ${\rm Sol}{\rm (VP)}=S$. 

Observe that ${\rm Sol}{\rm (VP)}$ is the union of two semi-closed generalized polyhedral convex sets. Furthermore, ${\rm Sol}{\rm (VP)}$ is \textit{disconnected} and \textit{non-closed}.

\medskip\noindent
{\sc Claim 2.} \textit{It holds that}	
\begin{equation}\label{ex_wsol_nonconvex}
\begin{aligned}
{\rm Sol}^w{\rm (VP)}&=\left\{u \in L \mid \langle x^*_1, u \rangle \leq 0, \, \langle x^*_2, u \rangle = 1\right\}\\
&\hspace*{2cm} \cup \left\{u \in L \mid \langle x^*_1, u \rangle = 0, \, \langle x^*_2, u \rangle \leq -1  \right\}.
\end{aligned}
\end{equation}

\smallskip
First, to clear that ${\rm Sol}^w{\rm(VP)} \subset S^w$, where $S^w$ is the set on the right-hand side of \eqref{ex_wsol_nonconvex}, take any $u \in {\rm Sol}^w{\rm (VP)}$. Suppose that $u=u_0+t_1e_1+t_2e_2$ with $u_0 \in X_0$ and $t_1, t_2 \in \mathbb{R}$. Since $u \in D$, $t_1 \leq 0$ and $t_2 \leq 1$. Observe that $u_0(t)=e_0(t)$ for all $t \in [-1,0]$; so, $u_0 \in D$. The inclusion $u \in {\rm Sol}^w{\rm (VP)}$ implies that $f(u)-f(x) \notin {\rm int}\,K$ for all $x \in D$. 

If $t_2 \geq 0$, then $u \in P_1$ and $f(u)=u+Tu=u_0+(t_1+t_2)e_1+t_2e_2$. Since the vector $x:=u_0+e_2$ belongs to $D \cap P_1$,  $f(x)=x+Tx=u_0+e_1+e_2$. Therefore, 
\begin{equation*}
f(u)-f(x)=(t_1+t_2-1)e_1+(t_2-1)e_2.
\end{equation*}
If $t_2 <1$, then $t_1+t_2 -1 < 0$ and $t_2-1 < 0$. So, $f(u)-f(x) \in{\rm int}\,K$ by~\eqref{ex_rep_intK}. This contradicts  the assumption $u \in {\rm Sol}^w{\rm (VP)}$. We thus get $t_2=1$. Consequently, $u \in S^w$.

If $t_2 < 0$, then $u \in P_2$ and $f(u)=u-Tu=u_0+(t_1-t_2)e_1+t_2e_2$. Clearly, the vector $x:=u_0+e_2$ belongs to $D \cap P_1$. Since $f(x)=x+Tx=u_0+e_1+e_2$,  
\begin{equation*}
f(u)-f(x)=(t_1-t_2-1)e_1+(t_2-1)e_2.
\end{equation*}
If $t_1-t_2-1 < 0$, then $f(u)-f(x) \in{\rm int}\,K$ by \eqref{ex_rep_intK}. This contradicts the assumption $u \in {\rm Sol}^w{\rm (VP)}$. It follows that $t_1-t_2-1 \geq 0$. Hence, $t_2 \leq -1$ as $t_1 \leq 0$. Therefore, $u \in S^w$. 

We have proved that ${\rm Sol}^w{\rm(VP)} \subset S^w$. To obtain the opposite inclusion, take any $u \in S^w$. Let $u_0 \in X_0$ and $t_1, t_2 \in \mathbb{R}$ be such that $u=u_0+t_1e_1+t_2e_2$. Given any $x \in D$, one can find a vector $x_0 \in X_0$, numbers $\tau_1 \leq 0$ and $\tau_2 \leq 1$ such that $x=x_0+\tau_1e_1+\tau_2e_2$. 

If $t_2=1$ and $t_1 \leq 0$, then $f(u)=u+Tu=u_0+(t_1+1)e_1+e_2$ by $u \in P_1$. If $0 \leq \tau_2 \leq 1$, then $x \in P_1$ and $f(x)=x_0+(\tau_1+\tau_2)e_1+\tau_2e_2$. Since 
\begin{equation*}
f(u)-f(x)=(u_0-x_0)+(t_1+1-\tau_1-\tau_2)e_1+(1-\tau_2)e_2
\end{equation*}
with $1-\tau_2 \geq 0$, one gets $f(u)-f(x) \notin{\rm int}\,K$ by \eqref{ex_rep_intK}. If $\tau_2 <0$, then $x \in P_2$ and $f(x)=x_0+(\tau_1-\tau_2)e_1+\tau_2e_2$. In accordance with \eqref{ex_rep_intK}, since 
\begin{equation*}
f(u)-f(x)=(u_0-x_0)+(t_1+1-\tau_1+\tau_2)e_1+(1-\tau_2)e_2
\end{equation*}
with $1-\tau_2 \geq 0$, one gets $f(u)-f(x) \notin{\rm int}\,K$. It follows that $f(u)-f(x) \notin{\rm int}\,K$ for all $x \in D$. Hence, $u \in {\rm Sol}^w{\rm (VP)}$. 

If $t_2 \leq -1$ and $t_1=0$, then  $f(u)=u-Tu=u_0-t_2e_1+t_2e_2$ by $u \in P_2$. If $0 \leq \tau_2 \leq 1$, then $x \in P_1$ and $f(x)=x_0+(\tau_1+\tau_2)e_1+\tau_2e_2$. Therefore,
\begin{equation*}
f(u)-f(x)=(u_0-x_0)+(-t_2-\tau_1-\tau_2)e_1+(t_2-\tau_2)e_2
\end{equation*}
By \eqref{ex_rep_intK}, since $-t_2-\tau_1-\tau_2 \geq 1-\tau_1-\tau_2 \geq 0$, one has $f(u)-f(x) \notin{\rm int}\,K$. If $\tau_2 <0$, then $x \in P_2$ and $f(x)=x_0+(\tau_1-\tau_2)e_1+\tau_2e_2$. As
\begin{equation*}
f(u)-f(x)=(u_0-x_0)+(-t_2-\tau_1+\tau_2)e_1+(t_2-\tau_2)e_2
\end{equation*}
with $(-t_2-\tau_1+\tau_2)+(t_2-\tau_2)=-\tau_1 \geq 0$, by \eqref{ex_rep_intK}, $f(u)-f(x) \notin{\rm int}\,K$. It follows that $f(u)-f(x) \notin{\rm int}\,K$ for all $x \in D$. Hence, $u \in {\rm Sol}^w{\rm (VP)}$. 
We have proved that $ {\rm Sol}^w{\rm (VP)}=S^w$. 

Clearly, ${\rm Sol}^w{\rm (VP)}$ is \textit{disconnected} and it is the union of two generalized polyhedral convex sets. 
\end{example}

\begin{acknowledgements}
This research was supported by the Vietnam National Foundation for Science and Technology Development (NAFOSTED) under grant number 101.01-2014.37. The author would like to thank Professor Nguyen Dong Yen for his guidance and the anonymous referees for valuable suggestions.
\end{acknowledgements}

\end{document}